\documentclass[12pt]{article}
\usepackage{graphicx}
\usepackage{amssymb}
\usepackage{amsmath}
\numberwithin{equation}{section}

\begin{document}
\author{Lev Sakhnovich}
\date{May 15, 2006}
\textbf{Levy Processes, Generators}

\begin{center} Sakhnovich Lev \end{center}

E-mail address: lev.sakhnovich@verizon.net \\

{\bf Abstract.} For a broad class of the Levy processes the new
form (convolution type) of the infinitesimal generators is
introduced. It leads to the new notions: a truncated generator, a
quasi-potential. The probability of the Levy process remaining
within the given domain is estimated.

\section{Main notions}
 Let us consider the Levy processes $X_{t}$ on $R$. If
$P(X_{0}=0)=1$ then Levy-Khinchine formula gives [1],[6]
\begin{equation}
\mu(x,t)=E\{\mathrm{exp}[izX_{t}]\}= \mathrm{exp}[-t\lambda(z)],\quad
t{\geq}0, \end{equation} where \begin{equation}
\lambda(z)=\frac{1}{2}Az^{2}-i{\gamma}z-\int_{-\infty}^{\infty}(e^{ixz}-1-ixz1_{D(x)})d\nu(x).
\end{equation} Here $A{\geq}0,\quad \gamma=\overline{\gamma},$ and
$D=\{x:|x|{\leq}1\}$ is closed unit ball, $\nu(x)$ is
monotonically increasing function  satisfying the conditions
\begin{equation}
\int_{-\infty}^{\infty}\frac{x^{2}}{1+x^{2}}d\nu(x)<\infty.
\end{equation} By $P_{t}(x_{0},B)$ we denote the probability
$P(X_{t}{\in}B)$ when $P(X_{0}=x_{0})=1$ and $B{\in}R$. The
transition operator is defined by the formula \begin{equation}
P_{t}f(x)=\int_{-\infty}^{\infty}P_{t}(x,dy)f(y).\end{equation}
Let $C_{0}$ be the Banakh space of continuous functions $f(x)
\quad (-\infty<x<\infty)$ satisfying the condition
$\mathrm{lim}f(x)=0,\quad |x|{\to}\infty$  with the norm $||f||=
\mathrm{sup}_{x}|f(x)|$. We denote by $C_{0}^{n}$ the set of
$f(x){\in}C_{0}$ such that $f^{k}(x){\in}C_{0},\quad
(1{\leq}k{\leq}n).$ It is known that [6]
\begin{equation} P_{t}f{\in}C_{0}\end{equation}
if $f(x){\in}C_{0}.$\\
Now we formulate the following important result (see [2],[6]) .\\
\textbf{Theorem 1.1.}\emph{The family of the operators $P_{t}\quad
(t{\geq}0)$ defined by the Levy process $X_{t}$ is a strong
continuous  semigroup on $C_{0}$ with norm $||P_{t}||=1$. Let $L$
be its infinitesimal generator. Then}
\begin{equation} Lf=\frac{1}{2}A\frac {d^{2}f}{dx^{2}}+\gamma
\frac{df}{dx}+\int_{-\infty}^{\infty}(f(x+y)-f(x)-y\frac{df}{dx}1_{D(x)})d\nu(x),\end{equation}
\emph{where} $f{\in}C_{0}^{2}$.
\section{Convolution type form of infinitesimal generator}
1. In this section we prove that under some conditions the
infinitesimal generator $L$ can be represented in the special
convolution type form \begin{equation}
Lf=\frac{d}{dx}S\frac{d}{dx}f,\end{equation} where the operator
$S$ is defined by the relation \begin{equation}
Sf=\frac{1}{2}Af+\int_{-\infty}^{\infty}k(y-x)f(y)dy,\end{equation}
and for arbitrary $M (0<M<\infty)$ we have
\begin{equation} \int_{-M}^{M}|k(t)|dt<\infty.\end{equation}
The representation $L$ in form (2.1) is convenient as the operator
$L$ is expressed with the help of the classic differential and
convolution operators.\\By $C_{c}$ we denote the set of functions
$f(x){\in}C_{0}$ with compact support.\\ \textbf{Lemma 2.1.}
\emph{Let the following conditions be
fulfilled.\\
1. The function $\nu(x)$ is monotonically increasing, has the
derivative when $x{\ne}0$ and}
\begin{equation}
\int_{-\infty}^{\infty}\frac{x^{2}}{1+x^{2}}d\nu(x)=
\int_{-\infty}^{\infty}\frac{x^{2}}{1+x^{2}}{\nu}^{\prime}(x)dx<\infty,
\end{equation}
\begin{equation}\nu(x){\to}0, \quad x{\to}\infty.\end{equation}
\emph{2. For arbitrary $M  (0<M<\infty)$ we have}
\begin{equation}
\int_{-M}^{M}|\nu(x)|dx<\infty,\quad
\int_{-M}^{M}|x|{\nu}^{\prime}(x)dx<\infty.\end{equation} 3.
\begin{equation}
x\nu(x){\to}0,\quad x{\to}0.\end{equation} \emph{Then the
equality}
\begin{equation}J=
\int_{-\infty}^{\infty}[f(y+x)-f(x)]{\nu}^{\prime}(y)dy,
\end{equation}
\emph{is true, where} $f(x){\in}C_{0}^{2},\quad
k(x)=\int_{0}^{x}\nu(y)dy$ and
\begin{equation}
J=\frac{d}{dx}\int_{-\infty}^{\infty}f^{\prime}(y)k(y-x)dy.\end{equation}
\emph{Proof.} For every $f(x){\in}C_{c}$  there exists such $M
(0<M<\infty)$ that
\begin{equation} f(x)=0 ,\quad x{\notin}[-M,M].\end{equation}Let
us introduce the following notations \begin{equation}
J_{1}=\frac{d}{dx}\int_{-\infty}^{x}f^{\prime}(y)k(y-x)dy,\end{equation}
\begin{equation}
J_{2}=\frac{d}{dx}\int_{x}^{\infty}f^{\prime}(y)k(y-x)dy.\end{equation}
Using (2.11) we have
\begin{equation}J_{1}=
-\frac{d}{dx}\int_{-M}^{x}[f(y)-f(x)+f(x)]k^{\prime}(y-x)dy.\end{equation}
From (2.11) and (2.13) we deduce the relation \begin{equation}
J_{1}=f(x)k^{\prime}(-M-x)+\int_{-M}^{x}[f(y)-f(x)]k^{\prime\prime}(y-x)dy.\end{equation}
When $M{\to}\infty$ we obtain the equality
\begin{equation}
J_{1}=\int_{-\infty}^{0}[f(y+x)-f(x)]k^{\prime\prime}(y)dy.\end{equation}
In the same way we deduce the relation \begin{equation}
J_{2}=\int_{0}^{\infty}[f(y+x)-f(x)]k^{\prime\prime}(y)dy.
\end{equation}The relation (2.9) follows directly from (2.15) , (2.16) and
equality
 $J=J_{1}+J_{2}.$
 The lemma is proved.\\
\textbf{Lemma 2.2.} \emph{Let the following conditions be
fulfilled.\\1. The function $\nu(x)$ satisfies conditions of Lemma
2.1.\\
2. For arbitrary $M\quad(0<M<\infty)$ we have
\begin{equation}
\int_{-M}^{M}|k(x)|dx<\infty,\quad
\int_{-M}^{M}|x{\nu}(x)|dx<\infty,\end{equation}where}
\begin{equation}k^{\prime}(x)=\nu(x),\quad x{\ne}0.\end{equation}
 3.
\begin{equation}
xk(x){\to}0,\quad x{\to}0;\quad x^{2}\nu(x){\to}0,\quad
x{\to}0.\end{equation} \emph{Then the equality}
\begin{equation}
J=\int_{-\infty}^{\infty}[f(y+x)-f(x)-y\frac{df(x)}{dx}1_{D(y)}]{\nu}^{\prime}(y)dy
+{\Gamma}f^{\prime}(x),
\end{equation}
\emph{is true, where} $\Gamma=\overline{\Gamma}$ and
$f(x){\in}C_{c}$.\\
\emph{Proof.} From (2.11)  we obtain  the relation
\begin{equation}
J_{1}=f^{\prime}(x)\gamma_{1}-\int_{x-1}^{x}[f^{\prime}(y)-f^{\prime}(x)]k^{\prime}(y-x)dy-
\int_{-M}^{x-1}f^{\prime}(y)k^{\prime}(y-x)dy,\end{equation} where
$\gamma_{1}=k(-1).$ We introduce the notations \begin{equation}
P_{1}(x,y)=f(y)-f(x)-(y-x)f^{\prime}(x),\quad
P_{2}(x,y)=f(y)-f(x).\end{equation}Integrating by parts (2.21) and
passing to limit when $M{\to}\infty$ we deduce that
\begin{equation}J_{1}=f^{\prime}(x)\gamma_{2}+\int_{x-1}^{x}P_{1}(x,y)k^{\prime\prime}(y-x)dy
+\int_{-M}^{x-1}P_{2}(x,y)k^{\prime\prime}(y-x)dy,\end{equation}
where $\gamma_{2}=k(-1)-k^{\prime}(-1).$ It follows from (2.22) and (2.23)
that \begin{equation}
J_{1}=\int_{-\infty}^{x}[f(y+x)-f(x)-y\frac{df(x)}{dx}1_{D(y)}]{\nu}^{\prime}(y)dy
+{\gamma}_{2}f^{\prime}(x).\end{equation}In the same way can be
proved the relation \begin{equation}
J_{2}=\int_{x}^{\infty}[f(y+x)-f(x)-y\frac{df(x)}{dx}1_{D(y)}]{\nu}^{\prime}(y)dy
+{\gamma}_{3}f^{\prime}(x),\end{equation}where
$\gamma_{3}=-k(1)+k^{\prime}(1).$ The relation (2.20) follows directly
from  (2.24) and (2.25). Here $\Gamma=\gamma_{2}+\gamma_{3}.$ The lemma is
proved.\\
\textbf{Remark 2.1.} The operator $L_{0}f=\frac{df}{dx}$ can be
represented in  form (2.1), where \begin{equation}
S_{0}f=\int_{-\infty}^{\infty}p_{0}(y-x)f(y)dy,\end{equation}
\begin{equation}p_{0}(x)=\frac{1}{2}\,\mathrm{sign}(x).\end{equation}\\
From Lemmas 2.1 , 2.2 and Remark 2.1 we deduce the following
assertion.\\
\textbf{Theorem 2.1.} \emph{Let the conditions of Lemma 2.1 or
Lemma 2.2 be fulfilled . Then the corresponding operator $L$ has
form }(2.1),(2.2).\\
\textbf{Proposition 2.1.} \emph{The generator $L$ of the Levy
process $X_{t}$admits the convolution representation $(2.1),(2.2)$ if
there exist such $C>0$ and $0<\alpha<2,\quad \alpha{\ne}1$ that}
\begin{equation}\nu^{\prime}(y){\leq}C|y|^{-\alpha-1},\quad
\end{equation} \emph{Proof.} The function $\nu(y)$  has the form
\begin{equation}
\nu(y)=\int_{-\infty}^{y}{\nu}^{\prime}(t)dt1_{y<0}-\int_{y}^{\infty}{\nu}^{\prime}(t)dt1_{y>0}.
\end{equation}We suppose that $1<\alpha<2$ and introduce the function
\begin{equation}
k_{0}(y)=\int_{-\infty}^{y}(y-t){\nu}^{\prime}(t)dt1_{y<0}-\int_{y}^{\infty}(t-y){\nu}^{\prime}(t)dt1_{y>0}.
\end{equation}We obtain the relation
\begin{equation}
k(y)=k_{0}(y)+(\gamma-{\Gamma})\mathrm{p}_{0}(y),\end{equation}
where $k_{0}(y)$ and $p_{0}(y)$ are defined by (2.27) and (2.30)
respectively.
 The constant $\Gamma$ is defined by relations:
\begin{equation}\Gamma=k_{0}(-1)-k_{0}^{\prime}(-1)-k_{0}(-1)+k_{0}^{\prime}(1),\quad 1{\leq}\alpha<2,
\end{equation}It follows from (2.28)-(2.30) that the conditions of Theorem
2.1 are fulfilled. Hence the proposition is true when
$1<\alpha<2.$ Let us consider the case when $0<\alpha<1$. In this
case we have
\begin{equation}
k_{0}(y)=\int_{-\infty}^{y}\nu^{\prime}(t)dt\,y+\int_{y}^{0}\nu^{\prime}(t)dt],
\quad y<0,\end{equation}
\begin{equation}
k_{0}(y)=-\int_{y}^{\infty}\nu^{\prime}(t)dt\,y-\int_{0}^{y}\nu^{\prime}(t)dt],
\quad y>0,\end{equation} and  \begin{equation}
k(y)=k_{0}(y)+\gamma\mathrm{p}_{0}(y)\quad 0<\alpha<1.\end{equation}In view of
(2.28) and (2.33),(2.34) the conditions of Theorem
2.1 are fulfilled. Hence the proposition is proved.\\
\textbf{Corollary 2.1.} \emph{If condition (2.28) is fulfilled then}
\begin{equation}
k_{0}(y){\geq}0,\quad
-\infty<y<\infty, \quad 1<\alpha<2,\end{equation}
\begin{equation}
k_{0}(y){\leq}0,\quad
-\infty<y<\infty, \quad 0<\alpha<1.\end{equation}
\textbf{Example 2.1.}The stable processes.\\
For the stable processes we have $A=0,\quad
\gamma=\overline{\gamma}$ and
\begin{equation}
\nu^{\prime}(y)=|y|^{-\alpha-1}(C_{1}1_{y<0}+C_{2}1_{y>0}),
\end{equation}where $C_{1}>0\quad C_{2}>0.$ Hence the function $\nu(y)$ has the form
\begin{equation}\nu(y)=\frac{1}{\alpha}|y|^{-\alpha}(C_{1}1_{y<0}-C_{2}1_{y>0}).
\end{equation}Let us introduce the functions
\begin{equation}k_{0}(y)=\frac{1}{\alpha(\alpha-1)}|y|^{1-\alpha}(C_{1}1_{y<0}+C_{2}1_{y>0}),
\end{equation} where $0<\alpha<2,\quad \alpha{\ne}1.$ When $\alpha
=1$ we have \begin{equation}k_{0}(y)=
-\mathrm{log}|y|\,(C_{1}1_{y<0}+C_{2}1_{y>0}).\end{equation} It
means that the conditions of proposition 2.1 are fulfilled. Hence
the generator $L$ for the stable processes admits the convolution
type
representation (2.1),(2.2).\\
\textbf{Proposition 2.2.} \emph{The kernel $k(y)$ of the operator
$S$ in representation $(2.1)$ for the stable processes has  form
$(2.31)$, when $1{\leq}\alpha<2$, and has form $(2.35)$ when $0<\alpha<1$.} \\
\textbf{Example 2.2.}The variance damped Levy processes [7].\\
For the variance damped Levy processes we have $A=0,\quad
\gamma=\overline{\gamma}$ and \begin{equation}
\nu^{\prime}(y)=C_{1}e^{-\lambda_{1}
|y|}|y|^{-\alpha-1}1_{y<0}+C_{2}e^{-\lambda_{2}|y|}y^{-\alpha-1}1_{y>0},
\end{equation}where
$C_{1}>0\quad C_{2}>0,\quad \lambda_{1}>0,\quad \lambda_{2}>0
\quad 0<\alpha<2,\quad \alpha{\ne}1.$ It follows from (2.42) that the
conditions of Proposition 2.1 are fulfilled.
Hence the generator $L$ for the variance damped Levy processes
admits the convolution type representation (2.1),(2.2) and the kernel $k(y)$
is defined by formulas (2.30),(2.31), when $1{\leq}\alpha<2$,and by formula
(2.35) when  $1{\leq}\alpha<2$.\\
\textbf{Example 2.3.}The variance Gamma process [7].\\
For the variance Gamma process we have $A=0,\quad
\gamma=\overline{\gamma}$ and
\begin{equation}
\nu^{\prime}(y)=C_{1}e^{-G
|y|}|y|^{-1}1_{y<0}+C_{2}e^{-M|y|}y^{-1}1_{y>0},
\end{equation}where $C_{1}>0\quad C_{2}>0,\quad G>0,\quad
M>0.$
 It follows
from (2.43) that the conditions of Proposition 2.1 are fulfilled and
the generator $L$ of variance Gamma process admits the convolution
type representation (2.1),(2.2).
 The kernel $k(y)$ is defined by relations
(2.33) and (2.34).\\
\textbf{Example 2.4.} The normal inverse Gaussian process [7].\\
In the case of the normal inverse Gaussian process we have
$A=0,\quad \gamma=\overline{\gamma}$ and
\begin{equation}
\nu^{\prime}(y)=Ce^{\beta y}K_{1}(|y|)|y|^{-1},\quad C>0,\quad
-1{\leq}\beta{\leq}1,\end{equation}where $K_{\lambda}(x)$ denotes
the modified Bessel function of the third kind with index
$\lambda$. Using equalities
\begin{equation}|K_{1}(|x|)|{\leq}Me^{-|x|}/|x|,\quad M>0,\quad
0<x_{0}{\leq}|x|,\end{equation}
\begin{equation}|K_{1}(|x|)x|{\leq}M,\quad
0{\leq}|x|{\leq}x_{0}\end{equation} we see that the conditions of
Proposition 2.1 are fulfilled. Hence the corresponding generator
$L$  admits the convolution type representation (2.1, {2.2}) and
 the kernel $k(y)$
is defined by  relations
(2.33) and (2.34).\\
\textbf{Example 2.5.} The Meixner process [7].\\
For the Meixner process we have \begin{equation}
\nu^{\prime}(y)=C\frac{\mathrm{exp}{\beta}x}{x\,\mathrm{sinh}{\pi}x}
,\end{equation}where $C>0,\quad -\pi<\beta<\pi.$ The conditions of
Proposition 2.1 are fulfilled. Hence the corresponding generator
$L$  admits the convolution type
representation (2.1),(2.2) and the kernel $k(y)$ is defined by relations
(2.33),(2.34).\\
\textbf{Remark 2.1.} Examples 2.1-2.5 are used in the finance
problems [7].\\
\textbf{Example 2.6.}Compound Poisson process.\\
We consider the case when $A=0, \quad \gamma=0$ and
\begin{equation}
M=\int_{-\infty}^{\infty}\nu^{\prime}(y)dy<\infty.\end{equation}
 Using formulas (2.1) and (2.2) we deduce that the corresponding
 generator
$L$ has the following form \begin{equation}
Lf=-Mf(x)+\int_{-\infty}^{\infty}\nu^{\prime}(y-x)f(y)dy.\end{equation}
\section{Potential}
The operator \begin{equation}
Qf=\int_{0}^{\infty}(P_{t}f)dt \end{equation}
is called \emph{potential} of the semigroup $P_{t}$. The generator $L$ and
the potential $Q$ are ( in general) unbounded operators. Therefore the operators
$L$ and $Q$ are defined not in the whole space $L^{2}(-\infty.\infty)$ but only
in the subsets $D_{L}$ and $D_{Q}$ respectively. We use the following property
of potential $Q$ (see[6]).\\
\textbf{Proposition 3.1.} \emph{If $f=Q\,g,\quad (g{\in}D_{q})$ then $f{\in}D_{L}$ and }
\begin{equation}
-Lf=g.\end{equation}
\textbf{Example 3.1.} Compound Poisson process.\\
Let the generator $L$ has form (2.49) where \begin{equation}
M=\int_{-\infty}^{\infty}\nu^{\prime}(x)dx<\infty,\quad \int_{-\infty}^{\infty}[\nu^{\prime}(x)]^{2}dx<\infty.
\end{equation} We introduce the functions \begin{equation}
K(u)=-\frac{1}{M\sqrt{2\pi}}\int_{-\infty}^{\infty}\nu^{\prime}(x)e^{-iux}dx,\end{equation}
\begin{equation}N(u)=\frac{K(u)}{1-\sqrt{2\pi}K(u)}.\end{equation}Let us note that
\begin{equation}
|K(u)|<\frac{1}{\sqrt{2\pi}},\quad u{\ne}0;\quad K(0)=-\frac{1}{\sqrt{2\pi}}.\end{equation}
It means that $N(u){\in}L^{2}(-\infty,\infty).$ Hence the function
\begin{equation}n(x)=-\frac{1}{\sqrt{2\pi}}\int_{-\infty}^{\infty}N(u)e^{-iux}du\end{equation}
belongs to $L^{2}(-\infty,\infty)$ as well. It follows from (2.49),(3.2)and (3.7) that the corresponding
potential $Q$ has form (see [8], Ch.11)\begin{equation}
Qf=\frac{1}{M}[f(x)+\int_{-\infty}^{\infty}f(y)n(x-y)dy].\end{equation}
\textbf{Proposition 3.2.} \emph{Let  conditions (3.3)  be fulfilled. Then
the operators $L$ and $Q$ are bounded in the space}  $L^{2}(-\infty,\infty).$\\
Now we shall give an example when the kernel $n(x)$ can be written in an explicit form.\\
\textbf{Example 3.2.} We consider the case when \begin{equation}
{\nu}^{\prime}(x)=e^{-|x|},\quad -\infty<x<\infty.\end{equation}
In this case $M=2$ and the operator $L$ takes the form
\begin{equation}
Lf=-2f(x)+\int_{-\infty}^{\infty}f(y)e^{-|x-y|}dy.\end{equation}
Formulas (3.4)-(3.7) imply that
\begin{equation}
Qf=\frac{1}{2}f(x)-\frac{1}{4\sqrt{2}}\int_{-\infty}^{\infty}f(y)e^{-|x-y|\sqrt{2}}dy.\end{equation}
\section{Truncated generators and quasi-potentials}
Let us denote by $\Delta$ the set of segments $[a_{k},b_{k}]$ such that\\
 $a_{1}<b_{1}<a_{2}<b_{2}<...<a_{n}<b_{n},\quad 1{\leq}k{\leq}n.$ By $C_{\Delta}$
 we denote the set of functions $g(x)$ on $L^{2}(\Delta)$ such that
\begin{equation}g(a_{k})=g(b_{k})=g^{\prime}(a_{k})=g^{\prime}(b_{k}),\quad 1{\leq}k{\leq}n,
\quad g^{\prime\prime}(x){\in}L^{2}(\Delta).
\end{equation}We introduce the operator $P_{\Delta}$ by relation
$P_{\Delta}f(x)=f(x)$ if $x{\in}\Delta$ and $P_{\Delta}f(x)=0$ if $x{\notin}\Delta.$\\
\textbf{Definition 4.1.} The operator
\begin{equation}
L_{\Delta}=P_{\Delta}LP_{\Delta} \end{equation}
is called \emph{a truncated generator}.\\
\textbf{Definition 4.2.} The operator $B$ with dense in $L^{2}(\Delta)$ definition domain $L_{B}$
is called \emph{a quasi-potential} if the functions $f=Bg,\quad (g{\in}L_{B})$ belong to $C_{\Delta}$
and \begin{equation}
-L_{\Delta}f=g.\end{equation}
It follows from (4.3) that
\begin{equation} -P_{\Delta}Lf=g,\quad ( f=Bg).\end{equation}
\textbf{Remark 4.1.} In a number of cases (see the next section) we need relation (4.4).
In these cases we can use the quasi-potential $B$, which is often simpler
than the corresponding potential $Q$.\\
\textbf{Remark 4.2.} The operators of type (4.2) are investigated in book ([5],Ch.2).\\
\textbf{Definition 4.3.} We call the operator $B$ a \emph{regular} one if the following
conditions are fulfilled.\\
1). The operator $B$ has the form \begin{equation}
Bf=\int_{\Delta}\Phi(x,y)f(y)dy,\quad f(y){\in}L^{2}(\Delta),\end{equation}
where the function $\Phi(x,y)$ can have a discontinuity only when $x=y$ and
\begin{equation}
|\Phi(x,y)|{\leq}M|x-y|^{-\beta},\quad 0{\leq}\beta<1.\end{equation}
2). The range of $B$ is dense in $L^{2}(\Delta)$.\\
It follows from condition 2) and relation (4.4) that the range of $L$ is dense
in $L^{2}(\Delta)$.
Further we assume that the quasi-potential $B$ is regular.

\section{The Probability of the Levy process remaining within the given domain}
In many theoretical and applied problems it is important to estimate the quantity
\begin{equation}p(t,\Delta)=P\{X_{\tau}{\in}\Delta\},\quad 0{\leq}\tau{\leq}t,\end{equation}
i.e. the the probability that a sample of the process $X_{\tau}$ remains inside $\Delta$
for $0{\leq}\tau{\leq}t.$ \\
To derive the integro-differential equations corresponding to Levy processes we
use the argumentation by Kac [3]  and our own argumentation (see [5]). Now we get rid
 of the requirement for the process to be stable.\\
 We suppose that \begin{equation}
 M(t)=\frac{1}{2\pi}\int_{-\infty}^{\infty}|\mu(x,t)|dx<\infty,\quad t>0\end{equation}
 and \begin{equation}
 \int_{0}^{1}M(t)dt<\infty.\end{equation}
 According to (1.2) we have the inequality \begin{equation}
 \mathrm{Re}[\lambda(x)]{\geq}0.\end{equation}
 \textbf{Remark 5.1.} It follows from (5.2) and (5.4) that \begin{equation}
\int_{-\infty}^{\infty}|\mu(x,t)|^{p}dx<\infty,\quad t>0,\quad p{\geq}1.\end{equation}
\textbf{Proposition 5.1.} ( see [], Ch.5.) \emph{Let  condition (5.2) be fulfilled. Then the
corresponding Levy process has the continuous density}\begin{equation}
\rho(x,t)=\frac{1}{2\pi}\int_{-\infty}^{\infty}e^{-ixz}\mu(z,t)dz,\quad t>0\end{equation}
\emph{and}
\begin{equation}\rho(x,t){\leq}M(t).\end{equation}
Now we introduce the sequence of functions
\begin{equation}
Q_{n+1}(x,t)=\int_{0}^{t}\int_{-\infty}^{\infty}Q_{0}(x-\xi,t-\tau)V(\xi)Q_{n}(\xi,\tau)d{\xi}d\tau,
\end{equation} where \begin{equation}Q_{0}(x,t)=\rho(x,t),\quad V(x)=1-1_{\Delta}.\end{equation}
For Levy processes the following relation \begin{equation}
Q_{0}(x,t)=\int_{-\infty}^{\infty}Q_{0}(x-\xi,t-\tau)Q_{0}(\xi,\tau)d{\xi}\end{equation} is true.
Using (5.8) and (5.10) we have \begin{equation}
0{\leq}Q_{n}(x,t){\leq}t^{n}Q_{0}(x,t)/n!.\end{equation}
Hence the series \begin{equation}
Q(x,t,u)=\sum_{n-0}^{\infty}(-1)^{n}u^{n}Q_{n}(x,t) \end{equation} converges.
The probabilistic meaning of $Q(x,t,u)$ is defined by the relation (see [4],Ch.4)
\begin{equation}E\{\mathrm{exp}[-u\int_{0}^{t}V(X(\tau))d\tau],c_{1}<X(t)<c_{2}\}=
\int_{c_{1}}^{c_{2}}Q(x,t,u)dx.\end{equation}The inequality $V(x){\geq}0$ and
relation (5.13) imply that the function $Q(x,t,u)$ monotonically decreases with respect
 to the variable "u" and the formulas \begin{equation}
 0{\leq}Q(x,t,u){\leq}Q(x,t,0)=Q_{0}(x,t)=\rho(x,t)\end{equation}
 are true. In view of (5.6), (5.7) and (5.14) the Laplace transform
\begin{equation}
\psi(x,s,u)=\int_{0}^{\infty}e^{-st}Q(x,t,u)dt,\quad s>0.\end{equation}
has the meaning.
According to (5.8) the function $Q(x,t,u)$ is the solution of the equation
\begin{equation}
Q(x,t,u)+u\int_{0}^{t}\int_{-\infty}^{\infty}\rho(x-\xi,t-\tau)V(\xi)Q(\xi,\tau,u)d{\xi}d\tau=
\rho(x,t)\end{equation}
Taking from both parts of (5.16) the Laplace transform and bearing in mind  (5.15) we obtain
\begin{equation}
\psi(x,s,u)+u\int_{-\infty}^{\infty}V(\xi)R(x-\xi,s)\psi(\xi,s,u)d\xi=R(x,s),
\end{equation} where
\begin{equation}
R(x,s)=\int_{0}^{\infty}e^{-st}\rho(x,t)dt.\end{equation}Multiplying both parts
(5.17) by $\mathrm{exp}(ixp)$ and integrating them with respect to
$x \quad (-\infty<x<\infty)$ we have \begin{equation}
\int_{-\infty}^{\infty}\psi(x,s,u)e^{ixp}[s+\lambda(p)+uV(x)]dx=1.\end{equation}
Here we use relations (1.1), (5.6) and (5.18). Now we introduce the function
\begin{equation}h(p)=\frac{1}{2\pi}\int_{\Delta}e^{-ixp}f(x)dx,\end{equation}
where the function $f(x)$ belongs to$C_{\Delta}$. Multiplying both parts  of (5.19)
by $h(p)$ and integrating them with respect to $p (-\infty<p<\infty)$ we deduce
the equality
\begin{equation}
\int_{-\infty}^{\infty}\int_{-\infty}^{\infty}\psi(x,s,u)e^{ixp}[s+\lambda(p)]h(p)dxdp=f(0).
\end{equation} We have used the relations \begin{equation}
V(x)f(x)=0,\quad -\infty<x<\infty,\end{equation}
\begin{equation}
\frac{1}{2\pi}\mathrm{lim}\int_{-N}^{N}\int_{\Delta}e^{-ixp}f(x)dxdp=f(0),\quad N{\to}\infty.
\end{equation}
Since the function $Q(x,t,u)$ monotonically decreases with respect to $"u"$ then
 by (5.15) this is true for the function $\psi(x,s,u)$ also. Hence there exists the limit
 \begin{equation}
 \psi(x,s)=\mathrm{\lim}\psi(x,s,u),\quad u{\to}\infty,\end{equation}where
 \begin{equation}
 \psi(x,s)=0,\quad x{\notin}\Delta.\end{equation}
 The probabilistic meaning of $\psi(x,s)$
 follows from the equality \begin{equation}
 \int_{-\infty}^{\infty}e^{-st}p(t,\Delta)dt=\int_{\Delta}\psi(x,s)dx.\end{equation}
 Using the properties of the Fourier transform and  conditions (5.24) , (5.25)
  we deduce from (5.21) the following assertion.\\
  \textbf{Proposition 5.2.} \emph{Let conditions (5.2) and (5.3)
  be fulfilled. Then the relation  \begin{equation}
  ((sI-L_{\Delta})f,\psi(x,s))_{\Delta}=f(0)\end{equation} is true.}\\
  The important function $\psi(x,s)$  can be expressed with the help of the
  quasi-potential $B$.\\
  \textbf{Theorem 5.1.} \emph{Let conditions (5.2), (5.3) be fulfilled and let the
  quasi-potential be regular.
  Then in the space $L(\Delta)$ there is one and only one function
  \begin{equation}
  \psi(x,s)=(I+sB^{\star})^{-1}\Phi(0,x) \end{equation}
  which satisfies relation (5.27).}\\
  \emph{Proof.} In view of (4.4) we have   \begin{equation}
  -BL_{\Delta}f=f,\quad f{\in}C_{\Delta}. \end{equation}Relations (5.28) and (5.29)
  imply that \begin{equation}
((sI-L_{\Delta})f,\psi(x,s))_{\Delta}=-((I+sB)L_{\Delta}f,\psi)_{\Delta}=-(L_{\Delta}f,\Phi(0,x))_{\Delta}.\end{equation}
Since $\Phi(0,x)=B^{\star}\delta(x),$  ($\delta(x)$  is the Dirac function)
then according to (5.24) and (5.30)
relation (5.27) is true. Let us suppose that in $L(\Delta)$ there is another
function $\psi_{1}(x,s)$
satisfying (5.27). Then the equality \begin{equation}
((sI-L_{\Delta})f,\phi(x,s))_{\Delta}=0,\quad \phi=\psi-\psi_{1} \end{equation}
is valid. We write  relation (5.31) in the form
\begin{equation}
(L_{\Delta}f,(I+sB^{\star})\phi)_{\Delta}=0.\end{equation}
The range of $L_{\Delta}$ is dense in $L^{\infty}(\Delta).$ Hence in view of (5.32) we have
$\phi=0.$
This proves the theorem.
\section{The kernel of the quasi-potential}
In this section we shall investigate the properties of the kernel $\Phi(x,y)$
of the quasi-potential $B$.\\
\textbf{Proposition 6.1.} \emph{Let conditions (5.2), (5.3) be fulfilled and let
the quasi-potential $B$ be regular. Then the corresponding kernel $\Phi(x,y)$
is non-negative i.e.} \begin{equation}\Phi(x,y){\geq}0.\end{equation}
\emph{Proof.} In view of (5.14) and (5.15) we have $ \psi(x,s,u){\geq}0$. Relation
(5.24)implies that  $ \psi(x,s){\geq}0$. Now it follows from (5.28) that
\begin{equation}\Phi(0,x)=\psi(0,x){\geq}0.\end{equation}
Let us consider the domains $\Delta_{1}$ and $\Delta_{2}$ which are connected
by relation $\Delta_{2}=\Delta_{1}+\delta$. We denote the corresponding truncated generators , quasi-potentials
kernels by $L_{k}, B_{k}$ and $\Phi_{k}(x,y)$, (k=1,2).
We introduce the unitary operator \begin {equation}
Uf=f(x-\delta), \end{equation}
which maps the space $L^{2}(\Delta_{2})$ onto $L^{2}(\Delta_{1})$. At the beginning
we suppose that the conditions of Theorem 2.1 are fulfilled.
Using formulas (2.1) and (2.2)we deduce that \begin{equation}
L_{2}=U^{-1}L_{1}U.\end{equation}
Hence the equality
\begin{equation}
L_{2}=U^{-1}L_{1}U.\end{equation}
Hence the equality
\begin{equation}
B_{2}=U^{-1}B_{1}U\end{equation} is valid. The last equality can be written in
the terms of the kernels \begin{equation}
\Phi_{2}(x,y)=\Phi_{1}(x+\delta,y+\delta).\end{equation}According to (6.2) and (6.6)
we have \begin{equation}\Phi_{1}(x+\delta,y+\delta){\geq}0.\end{equation}
As $\delta$  is an arbitrary real number relation (6.1) follows directly from (6.7).
 We remark that an arbitrary generator operator $L$ can be approximated by the operators of form
(2.1)(see[],Ch.). Hence the proposition is proved.\\
In view of (4.1), (4.5) and relation $Bf{\in}C_{\Delta}$  the following assertion is true.\\
\textbf{Proposition 6.2.} \emph{Let the quasi-potential $B$ be regular. Then the equalities}
\begin{equation}
\Phi(a_{k},y)=\Phi(b_{k},y)=0 \quad 1{\leq}k{\leq}n \end{equation}
\emph{are valid}.
\section{Sectorial operators}
We introduce the following notions.\\
\textbf{Definition 7.1.} The bounded operator $B$ in the space $L^{2}(\Delta)$ is called
\emph{sectorial} if \begin{equation}
(Bf,f){\ne}0,\quad f{\ne}0 \end{equation}
and
\begin{equation}
-\frac{\pi}{2}\beta{\leq}\mathrm{arg}(Bf,f){\leq}\frac{\pi}{2}\beta,\quad 0{\leq}\beta{\leq}1.
\end{equation}
\textbf{Definition 7.2.} The sectorial operator $B$ is called \emph{strong sectorial} if
$\beta<1$.\\
It is easy to see that the following assertion is true.\\
\textbf{Proposition 7.1.} \emph{ Let the operator $B$ be sectorial.
Then the operator $(I+sB)^{-1}$
is bounded when} $s{\geq}0$.

\begin{center} \textbf{References} \end{center}
1.Bertoin J., Levy Processes, University Press, Cambridge, 1996.\\
2. Ito K.,On Stokhastic Differential Equations, Memoirs Amer. Math. Soc. No.4,1951.\\
3.Kac M., On Some Connections Between Probability Theory and Differential and Integral
Equations,Proc.Second Berkeley Symp.Math.Stat. and Prob. Berkeley,189-215, 1951\\
4. Kac M., Probability and Related Topics in Physical Sciences, Colorado, 1957.\\
5.Sakhnovich L.A.,Integral Equations with Difference Kernels on Finite Intervals,
Operator Theory, v.84,Birkhauser,1996.\\
6.Sato K., Levy Processes and Infinitely Divisible Distributions, University
Press, Cambridge, 1999.\\
7. Schoutens W,Symens S., The Pricing of Exotic Options  by Monte-Carlo Simulation
in a Levy Market with Stokhastic Volatility, Preprint,1-27, 2002
8.Titchmarsh E.C., Introduction to the Theory of Fourier Integrals,Oxford, 1937
\end{document}